\documentclass[%
   , colorlinks % also like hyperref with CSC colorscheme (this is the default)
]{mpi2015-cscpreprint}

\usepackage[american]{babel}
\usepackage{amsthm}
\usepackage{cite}
\usepackage{amsmath,amssymb,amsfonts}
\usepackage{graphicx,bbm}

\usepackage{multicol}        % used for the two-column index
\usepackage[bottom]{footmisc}% places footnotes at page bottom

\usepackage{xcolor}
\usepackage{url}

\newcommand{\vertrule}[1][1.5ex]{\rule{0.4pt}{#1}}

\newcommand{\bA}{{\mathbf A}}

\newcommand{\bC}{{\mathbf C}}

\newcommand{\bM}{{\mathbf M}}

\newcommand{\bH}{{\mathbf H}}

\newcommand{\bX}{{\mathbf X}}
\newcommand{\bx}{{\mathbf x}}

\newcommand{\bV}{{\mathbf V}}

\newcommand{\cR}{ {\cal R} }

\newcommand{\SciML}{\textsf{SciML}}

\newcommand{\OpInf}{\textsf{OpInf}}
\newcommand{\SINDy}{\textsf{SINDy}}

\newcommand{\POD}{\textsf{POD}}

\newcommand*{\vertbar}{\rule[-1ex]{0.5pt}{2.5ex}}

\usepackage{algorithm}
\begin{document}

%%%%%%%%%%%%%%%%%%%%%%%%%%%%%%%%%%%%%%%%%%%%%%%%%%%%%%%%%%%%%%%%%%%%%%%%%%%%%%%%
% PAPER INFORMATION.                                                           %
%%%%%%%%%%%%%%%%%%%%%%%%%%%%%%%%%%%%%%%%%%%%%%%%%%%%%%%%%%%%%%%%%%%%%%%%%%%%%%%%

\title{Learning reduced-order Quadratic-Linear models in Process Engineering using Operator Inference}
  
\author[$\ast$ ]{Ion Victor Gosea}
\affil[$\ast$]{Max Planck Institute for Dynamics of Complex Technical Systems, Sandtorstr. 1, 39106 Magdeburg, Germany.\authorcr%
  \email{gosea@mpi-magdeburg.mpg.de}, \orcid{0000-0003-3580-4116}}% chktex 8
  
\author[$\dagger$ ]{Luisa Peterson}
\affil[$\dagger$]{Max Planck Institute for Dynamics of Complex Technical Systems, Sandtorstr. 1, 39106 Magdeburg, Germany.\authorcr%
  \email{peterson@mpi-magdeburg.mpg.de}, \orcid{0009-0006-0776-1517}}% chktex 8
  
\author[$\ddagger$]{Pawan Goyal}
\affil[$\ddagger$]{Max Planck Institute for Dynamics of Complex Technical Systems, Sandtorstr. 1, 39106 Magdeburg, Germany.\authorcr%
	\email{goyal@mpi-magdeburg.mpg.de}, \orcid{0000-0003-3072-7780}}%  

\author[$\ast \ast$]{Jens Bremer}
\affil[$\ast \ast$]{Clausthal University of Technology, Leibnizstr. 17, Clausthal-Zellerfeld, 38678, Germany.\authorcr%
	\email{bremer@icvt.tu-clausthal.de}, \orcid{0000-0002-9811-6698}}% chktex 8

\author[$\dagger \dagger$]{Kai Sundmacher}
\affil[$\dagger \dagger$]{Max Planck Institute for Dynamics of Complex Technical Systems, Sandtorstr. 1, 39106 Magdeburg, Germany.\authorcr%
	\email{sundmacher@mpi-magdeburg.mpg.de}, \orcid{0000-0003-3251-0593}}% chktex 8

\author[$\ddagger \ddagger$]{Peter Benner}
\affil[$\ddagger \ddagger$]{Max Planck Institute for Dynamics of Complex Technical Systems, Sandtorstr. 1, 39106 Magdeburg, Germany.\authorcr%
	\email{benner@mpi-magdeburg.mpg.de}, \orcid{0000-0003-3362-4103}}%  	

\shorttitle{Operator Inference for Process Engineering applications}
%\shortauthor{I. V. Gosea, L. Peterson, P. Goyal, J. Bremer, K. Sundmacher, P. Benner}
\shortauthor{I. V. Gosea et al.}
%\shortdate{}
  

%\msc{MSC1, MSC2, MSC3}
  
\abstract{%
  In this work, we address the challenge of efficiently modeling dynamical systems in process engineering. We use reduced-order model learning, specifically operator inference. This is a non-intrusive, data-driven method for learning dynamical systems from time-domain data. The application in our study is carbon dioxide methanation, an important reaction within the Power-to-X framework, to demonstrate its potential. The numerical results show the ability of the reduced-order models constructed with operator inference to provide a reduced yet accurate surrogate solution. This represents an important milestone towards the implementation of fast and reliable digital twin architectures.
  }

\novelty{This contribution illustrates the robustness and approximation capabilities of the operator inference method for a test case from process engineering.}

\maketitle

%%%%%%%%%%%%%%%%%%%%%%%%%%%%%%%%%%%%%%%%%%%%%%%%%%%%%%%%%%%%%%%%%%%%%%%%%%%%%%%%
% PAPER CONTENT.                                                               %
%%%%%%%%%%%%%%%%%%%%%%%%%%%%%%%%%%%%%%%%%%%%%%%%%%%%%%%%%%%%%%%%%%%%%%%%%%%%%%%%

\section{Introduction}\label{sec:Intro}

Dynamical models are particularly relevant for several important tasks and applications in the field of engineering science, such as analyzing transient behavior under operating conditions, enforcing parameter optimization, enabling long-time horizon prediction, and aiding control techniques. %For complex processes, obtaining a semi-discretized model (in space), describing all detailed dynamics and features is an important challenge.
Especially in the case of high-fidelity models, their state-space dimension can easily grow to the order of tens or even hundreds of thousands. In such cases, handling such models becomes computationally prohibitive with respect to storage or time constraints.
%Consequently, complexity reduction of the models is of great importance and need, and can facilitate relevant tasks such as repetitive simulations, controller design, and optimization.

Model order reduction (MOR) and reduced-order modeling (RoMod) refer to classes of methodologies that can be used to simplify complex dynamical models or to identify the underlying dynamics from data while enforcing an accurate approximation and using as little computational effort as possible. Traditionally, reduced-order models (ROMs) constructed through intrusive MOR approaches rely on the availability of the underlying mathematical equations (state-space realizations), i.e., on the actual models. For an overview of classical, intrusive (projection-based) approaches, we refer the reader to \cite{ACA05,BOCW17}. An advantage of such techniques is the existence of rigorous theoretical guarantees by means of a posteriori error estimation. A potential drawback is given by the fact that the ROMs computed via intrusive methods employ explicit projection of the governing equations onto low-dimensional dominant subspaces \cite{benner2015survey}. Hence, access to the underlying mathematical equations that produced the (high-fidelity) simulations is typically required. 

An alternative approach to classical (intrusive) methods of MOR, which rely on explicit access to a large-scale model, is the use of data-driven (non-intrusive) techniques based on measurements or simulation data. Unlike intrusive methods, data-driven RoMod does not require explicit knowledge of the model structure or matrices. Instead, low-order models can be directly constructed by using solely time-domain data, such as snapshots of the system's state-space evolution, along with snapshots of control inputs or of the observed outputs. Among such methods, Operator Inference (\OpInf) introduced in \cite{morPehW16}, has established itself in the last decade as a reliable RoMod and learning approach that exploits the inherent structure of physics-based models to capture the underlying system's dynamics efficiently. \OpInf \ constructs ROMs in a non-intrusive manner via a data-driven regression problem that learns reduced matrices from snapshot data. Typically, the learned dynamics of the ROM depend linearly and quadratically on the state variables (in reduced coordinates), although extensions to fitting other structures (cubic polynomial in \cite{khodabakhshi2022non} and non-polynomial in \cite{morBenGKetal20}), were also proposed. Even when the dynamics are characterized by such nonlinearities (high-degree polynomials or more generic, analytic ones), one can equivalently embed them into a quadratic manifold by means of lifting \cite{gu2011qlmor,morBenB15,morQiaKPetal20}.

The efficient modeling, simulation, and complexity reduction of dynamical systems in the field of process and chemical engineering has become pivotal in today's industrial landscape, especially with the ever-growing digitalization trend in modern plants and units in the age of Industry 4.0 \cite{morBreGFetal17}. Driven by this motivation, our work proposes an application of established reduced-order modeling techniques, focusing on the \OpInf \ method, for a test case involving a CO$_2$ methanation reactor \cite{zimmermann2022load}. The preliminary results reported here provide insights into the usability of the surrogate ROMs computed with \OpInf. The approximation quality attained offers promising prospects for further developments and extensions. Finally, this approach is well-suited for handling the intricacies of the highly nonlinear process under study. This is essential to ensure computational efficiency and approximation accuracy for applications in process engineering and other related disciplines.

This paper is organized as follows. After the introduction and motivation are set up, in Section \ref{sec:DDmeth} we provide a comprehensive account of state-of-the-art methods, together with their applications and extensions. Then, in Section \ref{sec:MethInvest}, we introduce the main method of interest, \OpInf, together with its theoretical background and some recent innovations. Section \ref{sec:App} illustrates the application of interest, a CO$_2$ methanation reactor characterized by coupled nonlinear PDEs. We present various numerical results when applying \OpInf \ to this test case, while in Section \ref{sec:Conc} the conclusions and future research endeavors are summarized.

\section{Data-driven methods} \label{sec:DDmeth}

Non-intrusive data-driven RoMod methods offer an alternative to conventional MOR methods, by means of constructing ROMs directly from data, bypassing the need for accessing an explicit realization. System identification (SI) methods \cite{Lju99} are considered the precursor of RoMod methods in the field of systems and control engineering. 
%By means of incorporating input-output data in the identification process, such methods were fairly prolific in the last decades. 
%The subspace identification method \cite{van2012subspace} is one of these and it estimates linear time-invariant state-space models using only samples of the input and output signals. Another class of methods is represented by Loewner-matrix approaches \cite{mayo2007fsg}, and use rational approximation techniques for identifying a ROM from frequency-domain data.
%Recently, SI methods that make use of full state access (snapshots of the state variables for different time steps) have emerged, including the family of 
Dynamic Mode Decomposition (DMD) \cite{kutz2016dynamic} is used to extract the dominant dynamic modes of a system and is closely connected to the concept of Koopman operator. %\cite{koopman1932dynamical}.
%a mathematical tool that describes the evolution of scalar observables (functionals of measurable state variables) in infinite time.
%Through the extended DMD approaches in \cite{williams2015data,brunton2016koopman}, it becomes possible to extend DMD to nonlinear systems, with the extra challenge of selecting appropriate observables such that the dynamics become close to linear.
In the field of Scientific Machine Learning (\SciML), several techniques have emerged to address the challenges and requirements of dealing with complex scientific problems, in the presence of simulation data. 
%Selecting the most appropriate \ML \ algorithm can be challenging due to the variety of algorithms, computer architectures, and \ML \ models available. 
Artificial Neural Network (ANN) architectures are particularly attractive for \SciML \ applications due to their ability to handle noisy or inexact data, perform automatic differentiation, and their flexibility as mesh-free models. Therefore, there is a growing interest in combining ANN-based solutions with traditional physics-based approaches to improve performance and reliability. Some recent works in this direction are \cite{raissi2019physics,lu2021learning,goyal2021lqresnet}. \SINDy, as introduced by \cite{brunton2016discovering}, utilizes scientific knowledge to improve the fitted model performance and constructs parsimonious models by selecting only a few terms from a library of candidate functions (sparse identification). 

%\textcolor{red}{Below, we provide a short account of \OpInf \ contributions.}

\OpInf \ involves fitting structured models in reduced coordinates by means of solving a regression problem formed from snapshot data (trajectories of the state evolution in the time domain). This approach provides a versatile framework for the analysis and modeling of continuous-time processes in chemical engineering. For example, the continuous-time behavior of the temperature and conversion profile of a chemical reactor can be accurately represented using \OpInf, as will be demonstrated in this paper. To the best of the authors' knowledge, the present paper, together with the recent papers \cite{LP23Escape} and \cite{mcquarrie2021data}, is one of the first to address the application of \OpInf \ to problems in process and chemical engineering. Specifically, \cite{LP23Escape} deals with CO$_2$ methanation, which is also the focus of this paper, while \cite{mcquarrie2021data} examines a single-injector combustion process in the context of aerospace engineering. \OpInf \ has also been successfully applied to test cases in computational fluid dynamics \cite{morBenGHetal22} and mechanics \cite{filanova2023operator}.

In recent years, \OpInf \ was extended also to cope with higher-order polynomial terms \cite{morQiaKPetal20}, or even with non-polynomial terms \cite{morBenGKetal20}. Additionally, parametric problems were treated \cite{yildiz2021learning,mcquarrie2023nonintrusive}. Then, contributions that propose regularizing the learning problem to enable performance on large-scale systems were proposed in \cite{mcquarrie2021data}, while problems with noisy or low-quality data sets were treated in \cite{uy2023active}. Dealing with differential-algebraic equations (DAEs) for incompressible flow problems was treated in \cite{morBenGHetal22}. The problem of stability of the ROMs constructed via \OpInf \ was tackled in \cite{goyal2023guaranteed,sawant2023physics,morPonGB24}.  In \cite{goyal2021lqresnet}, the authors propose combining the \OpInf \ approach with certain deep ANN approaches, to infer the unknown nonlinear dynamics of the system. For more references, detailed principles, and technical details of \OpInf, we refer the reader to the recent survey paper \cite{kramer2024learning}.

\section{The method under investigation: Operator Inference}  \label{sec:MethInvest}

\subsection{General description}
\label{sec:OpInf}

One of the strengths of the \OpInf \ methodology is that it makes use of the known physical structures and constraints at a qualitative level. That is, one can assume nonlinearities of quadratic structure (these appear naturally in many established flow problems, such as in the Navier-Stokes or Burgers' equations). 

Consider a quadratic full-order model (FOM) described by
\begin{equation}\label{eq:QuadMod}
	\dot \bx(t) =  \bA\bx(t) +  \bH(\bx(t)\otimes \bx(t)) + \bC,
\end{equation}
where $\bx(t) \in \mathbb{R}^n, \bA \in \mathbb{R}^{n \times n}, \bH \in \mathbb{R}^{n \times n^2}$, and $\bC \in \mathbb{R}^n$. Data preparation is an important step of the method, i.e., the matrix $\bX$ of snapshots of the state variables in (\ref{eq:Xdot_X}). As we will see below, we can also perform this step in reduced dimension. One needs to build, in the beginning, the following time-domain snapshot matrix for time instances $t_i>0$ and $\bx_i := \bx(t_i)$ with $1 \leq i \leq k$
\begin{align} \label{eq:Xdot_X}
	{\bX} & : =
	\left[
	\begin{array}{llll}
		~ \vertbar & ~\vertbar &        & ~\vertbar \\
		\bx_0      & \bx_1     & \cdots & \bx_k     \\
		~\vertbar  & ~\vertbar &        & ~\vertbar
	\end{array}
	\right] \in \mathbb{R}^{n \times (k+1)}.
\end{align}
Then we compute a projection matrix $\bV \in \mathbb{R}^{n \times r}$ using $r$ dominant Proper Orthogonal Decomposition (\POD) basis vectors (see Ch.3 in \cite{holmes2012turbulence} for more details) so that $\Vert \bX - \bV\bV^\top \bX \Vert \leq \tau$,  depending on a tolerance value $\tau>0$. Typically, the matrix $\bV$ is chosen by performing a Singular Value Decomposition (SVD) of $\bX$, and by assembling the leading $r$ left singular vectors of $\bX$.

By setting $\hat{\bx}(t):= \bV^\top  \bx(t)$, the intrusively calculated ROM via the \POD \ approach is computed via a Galerkin projection: $\hat{\bA} = \bV^\top \bA \bV, \hat{\bH} = \bV^\top \bH (\bV \otimes \bV)$, and $\hat{\bC} =  \bV^\top \bC$. However, the key feature of \OpInf \ is to infer these matrices in a non-intrusive manner, i.e., without using the FOM matrices in 
(\ref{eq:QuadMod}).

The reduced-order state data matrices can be put together using the compressed snapshots $\hat\bx_i = \bV^\top \bx_i $, $\hat\bx_i ^ \otimes = \hat\bx_i  \otimes \hat\bx_i$, and also the time-derivative data matrix can be put together from estimated snapshots $\dot{\hat\bx}_i$ using $\hat{\bX}$, e.g. by employing a time
derivative approximation, as follows:
\begin{align*}
	\hat{\bX} & : =
	\left[
	\begin{array}{llll}
		~ \vertbar & ~\vertbar &        & ~\vertbar \\
		\hat\bx_0  & \hat\bx_1 & \cdots & \hat\bx_k \\
		~\vertbar  & ~\vertbar &        & ~\vertbar
	\end{array}
	\right], \quad
	\hat{\bX}^\otimes  : =
	\left[
	\begin{array}{llll}
		~ \vertbar          & ~\vertbar          &        & ~\vertbar          \\
		\hat\bx_0 ^ \otimes & \hat\bx_1 ^\otimes & \cdots & \hat\bx_k^ \otimes \\
		~\vertbar           & ~\vertbar          &        & ~\vertbar
	\end{array}
	\right], \quad 	\dot{\hat{\bX}}  =
	\left[
	\begin{array}{llll}
		~ \vertrule      & ~\vertrule       &        & ~\vertrule       \\
		\hat{\dot\bx}_0 & \dot{\hat\bx}_1 & \cdots & \dot{\hat\bx}_k \\
		~\vertbar       & ~\vertbar       &        & ~\vertbar
	\end{array}
	\right].
\end{align*}

Then, a ROM of the form
\begin{equation*}
	\dot{\hat\bx}(t) =  \hat\bA\hat\bx(t) +  \hat\bH(\hat\bx(t)\otimes \hat\bx(t)) + \hat\bC,
\end{equation*}
can be obtained using projected data by solving the following optimization problem:
\begin{equation}\label{eq:LSprob}	\min_{\hat\bA,\hat\bH,\hat\bC}	\left\|\dot{\hat \bX} - \hat\bA \hat\bX - \hat\bH\hat\bX^\otimes - \hat\bC \right\|_{F}^2 + \alpha \cR(\hat\bA, \hat\bH, \hat\bC),
\end{equation}
where $\left\| \bM \right\|_{F}$ denotes the Frobenius (or Euclidean) norm of a matrix $\bM$, defined as the square root of the sum of the absolute squared entries of $\bM$.  Typically, the least squares (LS) problem in (\ref{eq:LSprob}) is ill-conditioned, and hence requires regularization techniques (such as the Tikhonov approach, see \cite{yildiz2021learning,mcquarrie2021data} in the context of \OpInf), i.e., encoded in the term  $\cR(\hat\bA, \hat\bH, \hat\bC)$ above, which is scaled by the regularization (or penalty) parameter $\alpha > 0$. We refer the reader to Section 4 in \cite{morBenGHetal22}, for details on implementation aspects for this "standard version" of \OpInf. 

\subsection{Solving the optimization problem}

It is to be noted that the FOM is never explicitly used, at any step, in the \OpInf \ algorithm. When $\alpha = 0$ (no regularization case), the LS problem can be written as
\begin{align}\label{eq:LSprob2}
	\min_{\hat\bA,\hat\bH,\hat\bC}	 
	\left\|\dot{\hat \bX} - \begin{bmatrix}
		\hat\bA & \hat\bH & \hat\bC
	\end{bmatrix} \mathcal{D} \right\|_{F}^2, \ \text{where} \ \mathcal{D} =  \begin{bmatrix} \hat{\bX} \\  \hat{\bX}^{\otimes} \\ \mathbf{1}^T
	\end{bmatrix},
\end{align}
with the matrix $\mathcal{D} \in  \mathcal{R}^{(r+r^2+1)\times (k+1)}$ constructed in terms of the projected states. Based on the singular value decay of $\mathcal{D}$ (and on a tolerance value), one can choose a suitable truncation order $\tilde{r}$. %different than the one indicated at the \POD \ step. 
%\textcolor{red}{Rephrase/Modify the following phrase!}
More precisely, one can write the solution of (\ref{eq:LSprob2}) as $\dot{\hat \bX} \mathcal{D}^{\dagger}$, where $\mathcal{D}^{\dagger}$ is the pseudo-inverse of matrix $\mathcal{D}$. This can be computed through a truncated SVD of $\mathcal{D}$, keeping only the $\tilde{r}$ dominant singular values and vectors  \cite{morBenGHetal22}. 

However, matrix  $\mathcal{D}$ may be ill-conditioned, in practice. Typically, a mixed SVD-regularization approach is to be used in practice, since choosing fewer singular values from $\mathcal{D}$ (small $\tilde{r}$) yields a larger mismatch of the data fidelity term, whereas, for larger $\tilde{r}$, the problem becomes more and more ill-conditioned. 
%An important property of \OpInf \ is that it recovers the intrusive \POD \ reduced model if data are Markovian, as shown in \cite{peherstorfer2020sampling}. 
A paramount innovation in the \OpInf \ framework was recently made in \cite{goyal2023guaranteed}, in which global-stability-informed learning was achieved through quadratic system parameterizations. There, the authors first provide ways of enforcing "local asymptotic stability", by imposing that the matrix $\hat{\bA}$ of the learned model has all its eigenvalues in the left-half complex plane (by imposing a parameterization of $\hat{\bA}$ in terms of symmetric positive definite (SPD) and skew-symmetric matrices). Then, an extension is proposed, that enforces "global asymptotic stability" by means of specially-tailored parameterization for the $\hat{\bH}$ operator of the learned quadratic model. Then,  a gradient-based approach is proposed in order to obtain a solution to the optimization problem in (\ref{eq:LSprob}), which now has additional constraints on the matrices involved. By relaxing the SPD condition, and by rewriting the problem of inferring operators via an integral form of the differential equation. This circumvents the need to incorporate derivative information, at the cost of a more complicated setup to be solved, involving integrals. By making use of a fourth-order Runge-Kutta scheme, the latter problem is approximated and solved iteratively, as described in detail in Section 6 of \cite{goyal2023guaranteed}.
In the numerics part, i.e., in Section \ref{sec:AppOpInf}, we leverage these results and provide additional explanations.

\section{Application under study}
\label{sec:App}

The application of interest, i.e., CO$_2$ methanation, represents an important process in the Power-to-X framework, which plays a key role in the storage of energy from renewable sources \cite{ghaib2018power}. CO$_2$ methanation, i.e., the conversion of carbon dioxide (CO$_2$) and hydrogen to methane and water, is a key aspect of this framework, facilitating the recycling of CO$_2$ emissions. The exothermic methanation reaction releases heat, which commonly leads to hotspots (localized areas of elevated temperatures) that affect temperature control and reactor efficiency. Effective management of these hotspots is critical, especially in reactors with variable energy inputs from renewable sources, such as wind or solar energy \cite{bremer2017co2}.

Our study investigates a one-dimensional, pseudo-homogeneous reactor model for catalytic CO$_2$ methanation, following  \cite{zimmermann2022load}. We employ a simplification by excluding the consideration of effective axial mass dispersion in our analysis. The behavior of the reactor is captured by coupled PDEs. These governing equations are written in terms of the energy balance, represented by the temperature variable $T$, and the mass balance, represented by the CO$_2$ conversion variable $X$, as:
%In our adaptation of the model, we assume a catalyst efficiency factor of 1, along with other modifications necessary for our specific analysis. The operational and design parameters of our model mirror those found in industrial reactors, enhancing the practical applicability of our research. These similarities allow our work to provide meaningful insights into CO$_2$ methanation as a viable approach for renewable energy storage.
%The governing equations for $X$ and $T$ with respect to axis $z$ over time $t$ are as:
\begin{align}\label{eq:PDEmodel}
	\varepsilon_\textrm{R} \frac{\partial X}{\partial t} &= -u \frac{\partial X}{\partial z} + \frac{M_{\mathrm{CO_2}}}{\rho y_{\mathrm{CO_2,in}}}(1-\varepsilon_\textrm{R}) \xi \sigma_{\textrm{eff}}, \\
	\left( \rho c_{\textrm{p}} \right)_{\textrm{eff}} \frac{\partial T}{\partial t} &= - u_{\textrm{in}} \rho_{\textrm{in}} c_\textrm{p} \frac{\partial T}{\partial z} + \frac{\partial }{\partial z} \left[ \Lambda_{\textrm{ax}} \frac{\partial T}{\partial z} \right] + \frac{4U}{D} (T - T_{\textrm{cool}}) - \Delta H_\textrm{R} (1-\varepsilon_\textrm{R}) \xi \sigma_{\textrm{eff}}, \nonumber
\end{align}
with initial and boundary conditions as given in \cite{zimmermann2022load}. In the formulation from (\ref{eq:PDEmodel}), various parameters appear such as the reaction rate $\sigma_{\textrm{eff}}$, the thermal gas conductivity $\Lambda_{\textrm{ax}}$, heat transfer coefficient $U$, etc. For more details, we refer the reader to \cite{zimmermann2022load}.

The mass and energy balances are coupled through various terms. For example, the reaction source term $\sigma_{\textrm{eff}}$ and the fluid density $\rho$ depend on both temperature and conversion. Additionally, fluid density $\rho$ influences the heat transfer coefficient $U$ and the mass flow $v_{\mathrm{gas}}$.

In our current study, we extend the analysis from our previous work \cite{LP23Escape} and concentrate on the start-up phase of the reactor.
Initially, the CO$_2$ conversion is zero throughout the reactor (denoted as $X(z,0) = 0$ for all $z$), and the internal temperature is uniformly set to the cooling temperature, $T_{\textrm{cool}}$. These initial conditions are crucial for analyzing the early dynamics of the reactor and provide a baseline for observing temporal changes. Formally, the initial conditions are
\begin{equation}
	X(z,0) = 0, \quad \forall z \in [0, L], \ \
	T(z,0) = T_{\textrm{cool}}, \quad \forall z \in [0, L],
\end{equation}
where $L$ is the length of the reactor. This setup allows a clear examination of the state evolution of the reactor from a consistent and well-defined starting point. In our previous study \cite{LP23Escape}, we investigated the dynamics of a reactor during a sudden increase in the cooling temperature, $T_{\textrm{cool}}$, from 270 to 280~°C.

Using a finite volume method, we semi-discretize the PDEs describing the reactor dynamics over $200$ equally-spaced control volumes. Such a high spatial resolution improves the fidelity of the dynamics captured but increases the computational requirements. In terms of temporal resolution, particular attention was paid to ensuring stability and convergence of the solution. The initial and boundary conditions align with the reactor start-up scenario mentioned earlier. Notably, alongside integrating at specific time points, we concurrently evaluate the right-hand side of the model to obtain accurate time derivatives. To solve the discretized equations, we used the \texttt{Kvaerno5} integrator from the \texttt{diffrax}~\footnote{See \url{https://github.com/patrick-kidger/diffrax} for details} library in \texttt{Python} 3.10, chosen for its balance between accuracy and computational efficiency. Computational performance, including run-time and resource utilization, was monitored to ensure the effectiveness of the simulation approach.

\begin{figure}[h!]
	\centering
	\includegraphics[width=0.8\columnwidth]{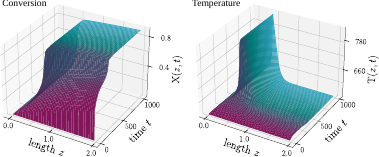}
	\caption{3D plots for CO$_2$ conversion (left) and temperature (right) over time and reactor length.}
	\label{fig:3D_plot}
\end{figure}

Figure~\ref{fig:3D_plot} depicts the evolution of conversion and temperature within the reactor during the start-up scenario, showcasing both the temporal and spatial domains. For the conversion variable $X$, the results show a progressive increase along the length of the reactor, reaching over 80~\% conversion towards the end. This indicates the effectiveness of the catalytic process in converting CO$_2$. The temperature variable $T$ varies between 550~K and just over 800~K during the simulation. A sharp peak in temperature is observed, representing the ignition of the reactor. This "shock-like" behavior is a critical feature of the start-up phase, demonstrating the dynamic nature of the system as it approaches a steady state. By the end of the simulation time, both the conversion and temperature profiles show convergence to a steady state, highlighting the ability of the reactor to stabilize under the given conditions.

\subsection{Numerical results}
\label{sec:AppOpInf}

\begin{figure}[h!]
	\centering
	\includegraphics[width=0.72\columnwidth]{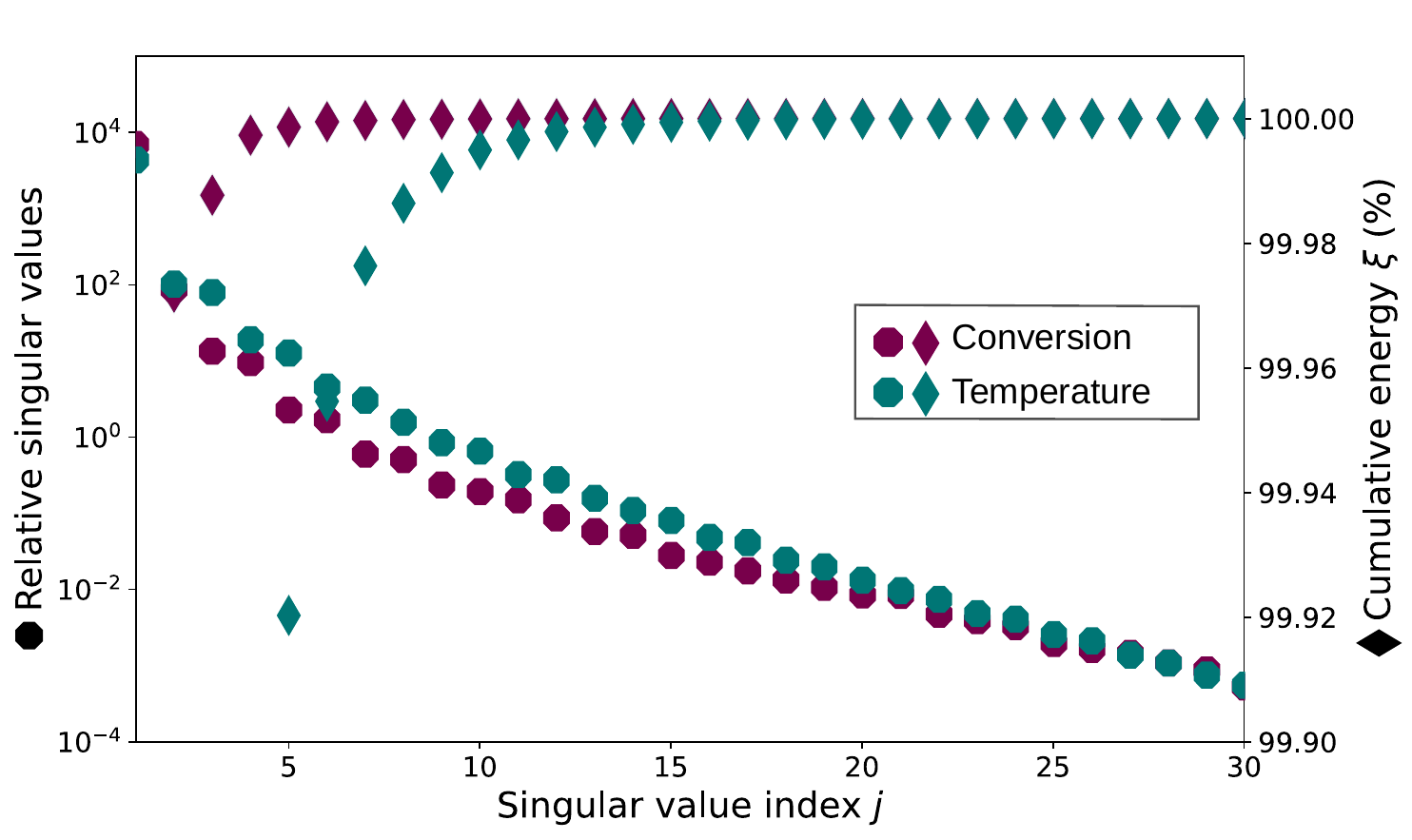}
	\caption{Representation of the singular value decay and the cumulative energy encompassed by the primary dominant modes.}
	\label{fig:SVD_results}
\end{figure}

In our study of the dynamic changes within the reactor, we applied Principal Component Analysis (PCA) via SVD to the collected state snapshots and initiated the \OpInf \ procedure. Figure~\ref{fig:SVD_results} shows the decay pattern of the singular values resulting from the PCA. Note that the decay of the singular values for the conversion variable $X$ is somewhat steeper than that for the temperature variable $T$. This indicates a more diverse distribution of singular values for the temperature, suggesting its more complex and non-linear characteristics compared to the conversion process. As a result, it is more difficult to capture the nuances of the temperature dynamics, reflecting its complexity. Aiming to capture 99.90~ to 99.99~\% of the total variance in the data, we developed ROMs across multiple dimensions. The variance of the data is measured by the cumulative sum of the squared singular values, which represents the information richness of the system state. The number of singular values required to reach this predefined variance threshold was determined by their cumulative contribution to the total variance.

In constructing our low-dimensional data and its derivatives, we projected the high-dimensional data onto the dominant modes using a reduced-order basis. This projection and subsequent operator inference were performed in \texttt{Python} 3.10, taking advantage of the language's computational efficiency and robust library support. Specifically, we used PyTorch's~\footnote{see \url{https://github.com/pytorch/pytorch} for details} library \texttt{Adam} as optimization algorithm combined with a \texttt{CyclicLR} scheduler for a cyclic learning rate policy. The scheduler, operating in a triangular2 mode with no cycle dynamics, adjusted the learning rate between $10^{-5}$ and $0.5$. The optimization was designed to stop early after 500 epochs if minimal improvement was observed, thereby increasing efficiency. A regularization factor of $10^{-4}$ was applied to the quadratic matrix for stability and to prevent over-fitting. The reduced operators  $\hat\bA$, $\hat\bH$ and $\hat\bC$ were derived following these methods, in accordance with the stability parameterization in \cite{goyal2023guaranteed}.

\begin{figure}[h!]
	\centering
	\includegraphics[width=0.72\columnwidth]{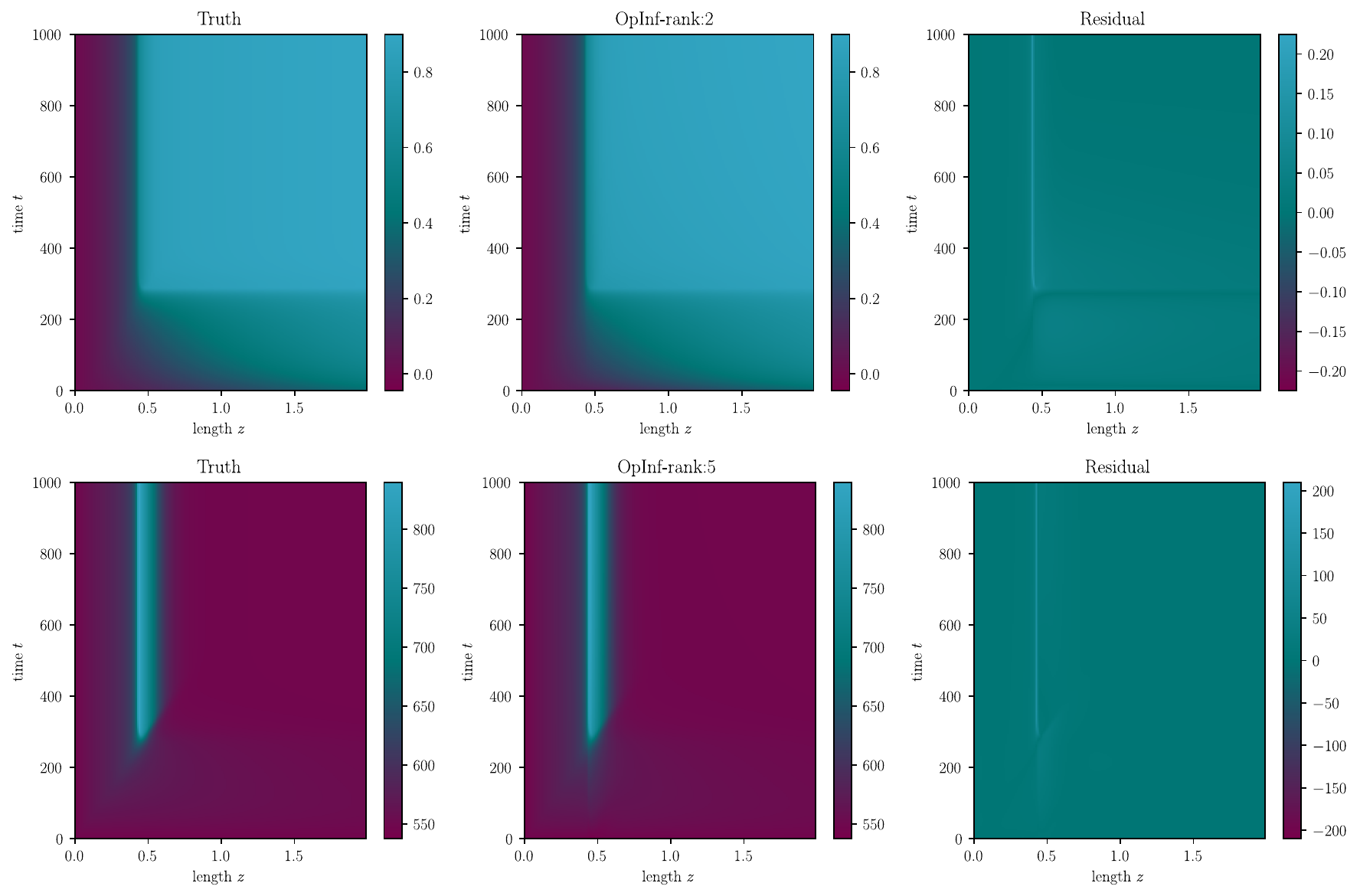}
	\caption{Flattened 3D representation comparing the true response with the inferred model, highlighting deviations in both conversion (first row) and temperature (second row).}
	\label{fig:flatted_3D}
\end{figure}

Figure~\ref{fig:flatted_3D} shows a flattened 3D plot comparing the actual data with our ROM, which captures 99.90~\% of the energy for both conversion and temperature. The model, at a rank of $r=7$ ($r_\mathrm{X}=2$ and $r_\mathrm{T}=5$), closely matches the actual data. However, discrepancies are observed in the regions representing the ignition phase, a near-shock scenario. Accurately capturing such abrupt and nonlinear transitions is a well-known challenge in modeling, often due to the inherent limitations of linear and quadratic terms in representing extreme state changes. These observations are critical for further refinement of the model, particularly in improving its ability to more accurately represent rapid dynamic changes. Quantitatively, the Frobenius norm shows a small deviation of 0.45~\% from the original model, underscoring the accuracy of the ROM. The computational efficiency of the ROM is high, in that it needs only 0.46~\% of the time required by the full mechanistic model.

\section{Conclusion and outlook}
\label{sec:Conc}

\OpInf \ combines the advantages of physics-based modeling with optimization and learning techniques by integrating first-principles models with data-driven regression.
The results reported in this work have successfully demonstrated the ability of the \OpInf \ ROMs to capture complex system dynamics in process engineering with increased computational efficiency. As an outlook, one option would be to allow variable parameters, by including variable input loads that affect species volume flow. In addition, the incorporation of specific control terms into the model will enhance its adaptability. By refining the tuning capabilities of the fitted models across a wide range of parameters, control settings, and many-query environments, they will hopefully be integrated into an effective digital twin environment.
%involving a small plant with many reactors in the near future.

\section*{Acknowledgement}

%\addcontentsline{toc}{section}{Acknowledgement}

This work is partially funded by the Bundesministerium für Bildung und Forschung (BMBF) and Project Management J\"ulich (PtJ) under grant 03HY302R. This work is also part of the research initiative “SmartProSys: Intelligent Process Systems for the Sustainable Production of Chemicals” funded by the Ministry for Science, Energy, Climate Protection and the Environment of the State of Saxony-Anhalt.

\vspace{3mm}

\addcontentsline{toc}{section}{References}
\bibliographystyle{spmpsci}
\bibliography{SciML_paper}

\end{document}